\documentclass[12pt,twoside,a4paper]{amsart}
\usepackage{amssymb}
\usepackage{amsthm}
\usepackage[english]{babel}
\usepackage[latin1]{inputenc} 
\usepackage[T1]{fontenc} 

\def\End{{\rm End}}

\def\deg{\mathrm{deg}\,}

\def\z{\zeta}

\def\w{\wedge}

\def\dbar{\bar\partial}
\def\ddbar{\partial\dbar}
\def\ddzi{\frac{\partial}{\partial \z _{i}}}

\def\C{{\mathbb C}}
\def\P{{\mathbb P}}

\def\Cn{\C^n}
\def\Pn{\P^n}

\def\D{{\mathcal D}}

\def\Hom{{\rm Hom\, }}

\def\End{\textrm{End}}

\def\be{\begin{equation}}
\def\ee{\end{equation}}

\newtheorem{thm}{Theorem}[section]
\newtheorem{lma}[thm]{Lemma}
\newtheorem{cor}[thm]{Corollary}
\newtheorem{prop}[thm]{Proposition}

\newtheorem{De}{Definition}

\newtheorem{preremark}{Remark}
\newtheorem{preex}{Example}

\begin{document}

\title{Explicit solutions of division problems for matrices of polynomials}
\author{Elin Götmark}

\maketitle

\begin{abstract}
We present integral representations of solutions to division problems in
$\mathbb{C}^n$ involving matrices of polynomials. We also find estimates
of the polynomial degree of the solutions by means of careful degree
estimates of the so-called Hefer forms which are components of the
representations. 
\end{abstract}

\tableofcontents

\section{Introduction}

Let $P(z) = (P_1(z), \ldots, P_m(z))$ be a tuple of polynomials in $\Cn$, and
$\Phi(z)$ another polynomial in $\Cn$ which vanishes on the zero set
$Z$ of $P$. By Hilbert's Nullstellensatz, we can
find $\nu \in \mathbb{N}$ and polynomials $Q(z) = (Q_1(z), \ldots,
Q_m(z))$ such that 
\[
\Phi^\nu = P \cdot Q = P_1 Q_1 + \cdots + P_m Q_m.
\]
There has been much research devoted to finding effective versions of this,
i.~e.~determining bounds on $\nu$ and the polynomial degrees of the 
$Q_i$'s. These degrees depend not only on
the degrees of $P$ and $\Phi$, but on $Z$ and the singularity of $P$ at
infinity. 

A breakthrough in the problem of degree estimates came in \cite{BR}, where
Brownawell 
found bounds by a combination of algebraic and analytical methods. 
The optimal degrees, which are slightly better than Brownawell's, were
found by Kollár in \cite{KO}, using purely algebraic 
methods. There are also classical theorems by M. Nöther and Macaulay
which treat simpler cases where one has conditions on the
singularity. 

In \cite{MA1}, the problem of solving 
\be \label{zabira}
\Phi = P \cdot Q = P_1 Q_1 + \cdots + P_m Q_m
\ee
with degree estimates is treated by means 
of residue currents on $\Pn$ based on the Koszul complex. One can then
recover the Nöther and Macaulay theorems. In the same article are explicit
solutions $Q$ constructed by means of integral representations; one then
has some loss of precision in the degree estimates. 

It is natural to look at generalizations of the division problem
(\ref{zabira}) to the case where $P$ is a matrix. More precisely, we
let $P$ be a polynomial 
mapping $\Cn \to \textrm{Hom} (\mathbb{C}^m, \mathbb{C}^r)$, i.~e.~an
$r \times m$ matrix $P$ where the entries are polynomials in $n$ complex 
variables. Let $\Phi = (\Phi_1, \ldots, \Phi_r)$ be an $r$-column of
polynomials. If we assume that there exists an $m$-column
of polynomials $Q$ such that $\Phi = P Q$, we want to find an explicit
solution $Q$ and get an estimate of its degree. 
We can reduce the case $r>1$ to the case $r=1$ by 
means of the Fuhrmann trick \cite{FU} (see Remark \ref{martha}),
however, we will then lose precision in the degree estimates.  

The case $r>1$ is treated in \cite{MA3}, where an estimate of the
degree of the solution is obtained by means of
residue currents on $\Pn$, based on the Buchsbaum-Rim complex. 
The aim of this paper is to present explicit integral representations
of $Q$ for $r \geq 1$. The paper \cite{MA2} 
contains a method of obtaining explicit solutions to division problems
$\Phi = P Q$ over open domains $X \subset \Cn$, and we will adapt this
method so that we can use it on $\Pn$. Also, we 
present careful degree estimates of so-called Hefer forms, which make it
possible to estimate the degree of the solutions. We will lose some
precision in the degree estimates compared to \cite{MA3}, but in
return we get explicit solutions. We also
point out that in principle, explicit solutions can be obtained which 
satisfy the same degree estimates as in \cite{MA3}.

\section{A division formula for tuples of polynomials} \label{tofsvipa}

Instead of solving the problem $\Phi = P Q$ in $\Cn$, we will
homogenize the expression and look at 
the corresponding problem in $\Pn$ instead. 
Let $\z = (\z_0, \z') = (\z_0, \z_1, \ldots \z_n)$. 
Set $\phi =\z_0^{\rho} \Phi(\z'/\z_0)$, where $\rho = \deg \Phi$; the
entries of $\phi$ will then be 
$\rho$-homogeneous holomorphic polynomials in $\C^{n+1}$, and can be
interpreted as sections of a line bundle over $\Pn$. Assume that 
$\deg(P^j) \leq d_j$, where $P^j$ is column number $j$ in $P$, and let 
$p(\z)$ be the matrix with columns
$p^j(\z) = \z_0^{d_j} P^j(\z'/\z_0)$. We will try to solve the division
problem $\phi = p \cdot q$ by means of integral
representation. Generally, this is not always possible, so we will get a
residue term $\phi R^p$ as well. When $\phi$ annihilates the
residue $R^p$, we can solve the division problem. Dehomogenizing $q$ will 
give a solution to our original problem. 

We will now find integral representation formulas which solve the
division problem in a more general setting, and return to our original
problem in Section \ref{xavier}. 
Let $S=\C[z_0, \ldots, z_n]$ be the ring of holomorphic homogeneous
polynomials over 
$\C^{n+1}$ with the ordinary grading, and let $S(-k)$ denote the same
ring, but where the grading is given by adding $k$ to the degree of the
polynomial. For example, the constants have degree $k$. Assume
that we have a graded complex $M$ of modules over $\mathbb{C}^{n+1}$ 
\[
0 \stackrel{f_{N+1}}{\longrightarrow} M_N \stackrel{f_N}{\longrightarrow} \cdots 
\stackrel{f_2}{\longrightarrow} M_1 \stackrel{f_1}{\longrightarrow} M_0 \to 0,
\]
where $M_m = S(-d_{m}^1) \oplus \ldots \oplus S(-d_{m}^{r_m})$ has
rank $r_m$ and basis $e_{m1}, \ldots,  e_{mr_m}$. We will take $d_0^i
= 0$ for all $i$. Since the complex is graded, $f_m$ must
preserve the degree of elements of $M_m$. To attain this, we must have 
\be \label{pinsent}
f_m = \sum_{i=1}^{r_{m-1}} \sum_{j=1}^{r_{m}} f_m^{ij} e_{m-1,i} \otimes
e_{mj}^\ast,
\ee
where $f_m^{ij}$ is a $d_{m-1}^i-d_m^j$-homogeneous holomorphic
polynomial. Let $M = M_0 + \cdots + M_N$. We define $f:M \to M$ to be
such that $f_{|M_m} = f_m$. 

Let $\mathcal{O}(-1)$ denote the tautological line bundle over $\Pn$, and if
$k$ is positive we let $\mathcal{O}(-k) : = \mathcal{O}(-1)^{\otimes k}$ and $\mathcal{O}(k) : =
(\mathcal{O}(-1)^\ast)^{\otimes k}$. Note that sections of $\mathcal{O}(k)$ can be seen
as $k$-homogeneous polynomials in $\C^{n+1}$. 
Our complex over $\C^{n+1}$ has an associated complex $E$ over $\Pn$,
namely 
\[
0 \to E_N \stackrel{f_N}{\longrightarrow} \cdots 
\stackrel{f_2}{\longrightarrow} E_1 \stackrel{f_1}{\longrightarrow} E_0 \to 0,
\]
where 
\[
E_m = \sum_{i = 1}^{r_m} \mathcal{O}(d_m^i) \otimes L_{mi}
\]
and the $L_{mi}$ are trivial line bundles with frames $e_{mi}$. We can
then write 
\[
f_m = \sum_{i=1}^{r_{m-1}} \sum_{j=1}^{r_{m}} f_m^{ij} e_{m-1,i} \otimes
e_{mj}^\ast,
\]
where $f_m^{ij}$ is the section of $\mathcal{O}(d_{m-1}^i-d_m^j)$ corresponding
to $f_m^{ij}$ in (\ref{pinsent}). Note that from $E$, we can obtain a
complex $E'$ of trivial 
vector bundles over $\Cn$, simply by taking the natural trivialization over
the coordinate neighborhood $\{\z_0 \neq 0\}$. We have 
\be \label{salieri}
0 \to E'_N \stackrel{F_N}{\longrightarrow} \cdots 
\stackrel{F_2}{\longrightarrow} E'_1 \stackrel{F_1}{\longrightarrow}
E'_0 \to 0, 
\ee
where $F_m(\z') = f_m(1,\z')$ for $\z' \in \Cn$, and $E'_m$ is a
trivial vector bundle of rank $r_m$ over $\Cn$. 

We introduce a $\mathbb{Z}_2$-grading on $E$, by writing  $E = E^+
\oplus E^-$ where $E^+ = \bigoplus E_{2k}$ and $E^- 
= \bigoplus E_{2k+1}$, and we say that sections of $E^+$ have even degree and
sections of $E^-$ have odd degree. This induces a grading
$\textrm{End} E = (\textrm{End} E)^+ \oplus (\textrm{End} E)^-$ where
a mapping $f:E \to E$ has odd degree if it maps $E^+ \to E^-$ and $E^-
\to E^+$, and similarly for even mappings. We also get an induced
grading on currents taking values in $E$, that is, $\D'_{\bullet,
  \bullet} (\Pn, E)$, by combining the gradings on
$\D'_{\bullet, \bullet} (\Pn)$ and $E$, so that $\alpha \otimes \phi =
(-1)^{\deg \alpha \cdot \deg \phi} \phi \otimes \alpha$ if $\alpha$ is
a current and $\phi$ a section of $E$. Similarly, we get
an induced grading on $\D'_{\bullet, \bullet} (\Pn, \End E)$, then $\dbar$
is an odd mapping. Since $f$ is odd and holomorphic, it is easy to
show that $\dbar \circ f = - f \circ \dbar$. 

Let $Z_i$ be the set where $f_i$ does not attain its highest rank, and set $Z
= \bigcup Z_i$. To construct integral formulas, we need to find
currents $U_k$ with bidegree $(0,k-1)$ and 
$R_k$ with bidegree $(0,k)$, with values in $\textrm{Hom}
(E_0,E_k)$, such that if $U = \sum_k U_k$ and $R = \sum_k R_k$, then 
\be \label{bone}
(f - \dbar)U = I_{E_0} - R. 
\ee
For the reader's convenience, we indicate the construction here, but
see \cite{AW} for details. We first define mappings
$\sigma_i:E_{i-1} \to E_i$ defined outside $Z_i$, such that $f_i
\sigma_i = \textrm{Id}$ on $\textrm{Im} f_i$, $\sigma_i$
vanishes on $(\textrm{Im} f_i)^\perp$, and $\textrm{Im} \sigma_i \perp
\textrm{Ker} f_i$. Now, set $u_1 = \sigma_1$ and 
$u_k = \sigma_k \dbar \sigma_{k-1} \w \cdots \w \dbar \sigma_1$. We have $(f -
\dbar)u = \sum_k (f_k u_k - \dbar u_k) = f_1 \sigma_1 + \sum_k
(\dbar u_{k-1} - \dbar u_k) = \textrm{Id}_{E_0}$ outside $Z$. 
To extend $u$ over $Z$, we let $U$ be the analytical continuation of
$|f|^{2\lambda} u$ to $\lambda = 0$, then we get 
(\ref{bone}). Note that for degree reasons we necessarily have $U_k
= 0$ for $k > \mu + 1$, where $\mu = \min(n,N-1)$. 

We will need to recall a proposition 
about representation of holomorphic sections of line bundles over
$\Pn$. See for example \cite{EG} for more on this. Let $z = (1,z')$ be
a fixed point, and set 
\[
\eta = 2 \pi i \sum_0^n z_i \ddzi.
\]
If we express a projective form in homogeneous coordinates, we can let
contraction with $\eta$ act on it and get a new projective form. In
fact, we will get an operator 
\[
\delta_\eta : \mathcal{D}'_{l+1,0}(\Pn,\mathcal{O}(k)) \to
  \mathcal{D}_{l,0}'(\Pn,\mathcal{O}(k-1)),
\]
where $\mathcal{D}'_{l,0}(\Pn,\mathcal{O}(k))$ denotes the
currents of bidegree $(l,0)$ taking values in $\mathcal{O}(k)$. 

\begin{prop} \label{hatsuyuki}
Let $\phi(\z)$ be a section of $\mathcal{O}(\rho)$, and let $g(\z,z) =
g_{0,0} + \ldots + g_{n,n}$ be a smooth form, where 
$g_{k,k}$ is of bidegree $(k,k)$ and takes values in
$\mathcal{O}(\rho + k)$. Assume that $(\delta_\eta - \dbar)g = 0$ 
  and $g_{0,0}(z,z) = \phi(z)$. Then we have 
\[
\phi(z) = \int_{\Pn}(g \w \alpha^{n+\rho})_{n,n},
\]
where
\[
\alpha = \alpha_0 + \alpha_1 = \frac{z \cdot \bar \z}{|\z|^2} - \dbar
\frac{\bar \z \cdot d\z}{2 \pi i |\z|^2}. 
\]
\end{prop}

The idea now is to find a weight $g$ that contains $f_1$ as a factor,
and then apply Proposition \ref{hatsuyuki}. As components of the
weight, we need so-called Hefer forms. 
We define $f^z_m:E_m \to E_{m-1}$ to be the mapping
\[
f^z_m = \sum_{i=1}^{r_{m-1}} \sum_{j=1}^{r_{m}} \z_0^{d_{m-1}^i-d_m^j}
f_m^{ij}(z) e_{m-1,i} \otimes e_{mj}^\ast.
\]

\begin{prop} \label{yak}
There exist $(k-l,0)$-form-valued mappings $h_k^l(\z,z): E_k \to E_l$,
such that $h_k^l = 0$ for $k<l$, $h_l^l = I_{E_l}$, and 
\[
\delta_\eta h_k^l = h_{k-1}^l f_k- f_{l+1}^z h_k^{l+1}
\]
Moreover, the coefficients of the form $(h_k^l)_{\alpha \beta}$, that
is, the coefficient of $e_{l\alpha} \otimes e_{k\beta}^\ast$, will take values in
$\mathcal{O}(d_{l}^\alpha-d_{k}^\beta+k-l)$, and can be chosen so that they are
holomorphic polynomials in $z'$ of degree $d_{l}^\alpha-d_{k}^\beta-(k-l)$. 
\end{prop}

We will prove this in the next section. Now we can define our weight:

\begin{prop} \label{eowyn}
If
\[
g = f_1^z \sum h_k^1 U_k + \sum h_k^0 R_k = f_1^z h^1 U + h^0 R,
\]
then $g_{0,0} {\arrowvert_{\Delta}} = I_{E_0}$ and $\nabla_\eta g = 0$. 
\end{prop}

The following proof is identical to the one in Section 5 in
\cite{MA2}, except that our proof will be in $\Pn$ instead of $\Cn$.
We include it for the reader's convenience. 
\begin{proof}
By definition, we have 
\be \label{signal}
\nabla_\eta g = f_1^z [(h^1 f - f^z h^2)U - h^1 \dbar U]
+ (h^0 f - f^{z} h^1)R - h^0 \dbar R.
\ee
Note that one has to check the total degree of $h^1$ and $h^0$ to get
all the signs correct. Since $\dbar U = fU-R-I_{E_0}$ and $f^2
= 0$, the right hand side of (\ref{signal}) is equal to
\[
f_1^z [h^1 fU - h^1(fU-R-I_{E_0})]+ h^0(fR-\dbar R)-f^{z} h^1R. 
\]
Furthermore, we have $h^1I_{E_0} = 0$ and $(f - \dbar)R = (f - \dbar)(R +
I_{E_0}) = (f - \dbar)^2 U = 0$, so it follows that $\nabla_\eta g = 0$. 
We also have $g_{0,0} = f^z h_1^1 U_1 = f^zU_1$, so that $g_{0,0}
{\arrowvert_{\Delta}} = fU_1 = I_{E_0}$. 
\end{proof}

This is the main result of this section:

\begin{thm} \label{inuvik}
Let $\phi(\z)$ be a section of $\mathcal{O}(\rho) \otimes E_0$. Fix $z =
(1,z')$ and let $\Phi(z') =
\phi(1,z')$, and $F_1(z') = f_1(1,z')$. If $g$ is defined as in
Proposition \ref{eowyn}, we have the following decomposition:
\[
\Phi(z') =  \int_\z g \phi \w \alpha^{n + \rho} = 
F_1(z') Q(z') + \int_\z h^0 R \phi \w \alpha^{n + \rho},
\]
where $Q(z') = \sum_1^{r_1} Q_i(z') e_{1j}$ with 
\[
Q_i(z') = \int_\z \z_0^{-d_1^i} (h^1 U \phi)_i(\z,z) \w \alpha^{n + \rho}.
\]
If $R \phi = 0$, the second integral will vanish and we
get a solution to our division problem. Moreover, we get the estimate
$\deg_{z'} F_1Q \leq \max_{k,\beta} (\rho - d_k^\beta)$.
\end{thm}

\begin{proof}
First, we note that 
\[
g \phi = \sum_{j=1}^{r_0} (g \phi)_j e_{0j},
\]
and according to Proposition \ref{eowyn}, $(g \phi)_j$ is
a weight of the type needed to apply Proposition \ref{hatsuyuki}. Now
we need only check the degree in $z'$. The term in $(g \phi)_j$ of
bidegree $(k,k)$ is $(f_1^z h_{k+1}^1 U_{k+1} \phi)_j$ (if we
disregard the terms containing $R$). We
must pair this with $\alpha_0^{\rho+k} \alpha_1^{n-k}$ to get
something with full bidegree in $z$. 

The degree in $z'$ of a term in 
\[
\sum_{k, \alpha, \gamma, \beta} 
(f_1^z)_{j \alpha} (h_{k+1}^1)_{\alpha \beta}(U_{k+1})_{\beta \gamma} \phi_\gamma \w
\alpha_0^{\rho +k} \alpha_1^{n-k}
\]
is $- d_1^\alpha + (d_1^\alpha - d_{k+1}^\beta - k) + (\rho 
+k) = \rho - d_k^\beta$. We get the estimate $\deg F_1 Q
\leq \max_{k,\beta} (\rho - d_{k+1}^\beta)$.
\end{proof}

\section{Constructing Hefer forms} \label{logan}

In this section we construct Hefer forms for the complexes $E$ and
$E'$, and investigate their degrees. We first state and prove a more
general theorem, and obtain as a corollary Hefer forms $H_k^l$ with
specific degrees for the complex $E'$ over $\Cn$. From the $H_k^l$ we
get corresponding Hefer forms $h_k^l$ for the complex $E$ over
$\Pn$, i.\ e.\ we prove Proposition \ref{yak} in the previous section. 

Let $\delta_{\z-z}$ denote contraction with the vector field  
\[
2 \pi i \sum_1^n (\z_j - z_j) \frac{\partial}{\partial \z_j}.
\]

\begin{thm} \label{linea}
Let $\phi(\z,z)$ be a $(l,0)$-form with holomorphic
polynomials of degree $r$ for coefficients. If $l>0$ 
and $\delta_{\z-z}\phi = 0$, we have $\phi = \delta_{\z-z} \psi$, where $\psi$
is a $(l+1,0)$-form with holomorphic polynomials of
degree $r-1$ for coefficients. If $l = 0$, let $\phi(\z)$ be a
holomorphic polynomial. We can then write $\phi(\z) - \phi(z) =
\delta_{\z-z} \psi$, where $\psi$ satisfies the same conditions as above. 
\end{thm}

For the proof of Theorem \ref{linea}, we need a lemma:

\begin{lma} \label{ineffable}
Let $\alpha = (\alpha_1, \ldots, \alpha_{n})$ and $\beta = (\beta_1,
\ldots, \beta_{n})$ be multiindices and $D_R$ the ball with
radius $R$ in $\Cn$. We have
\be \label{blodpudding}
I_{\alpha,\beta,R} := \int_{\partial D_R} w^\alpha \bar w^\beta
\partial|w|^2 \w ( \ddbar |w|^2 )^{n-1} = 0 
\ee
if $\alpha \neq \beta$ and (\ref{blodpudding}) is non-zero and proportional
to $R^{2(|\alpha|+n)}$ if $\alpha = \beta$.
\end{lma}

\begin{proof} (of Lemma \ref{ineffable})
The integral is unchanged if we take the pullback of the integrand by
the function $A(w) = (\lambda_1 w_1, \ldots
\lambda_n w_{n})$, where $\lambda_i \in \mathbb{C}$
and $|\lambda| = |(\lambda_1, \ldots \lambda_{n})| = 1$. We have
\[
I_{\alpha,\beta,R} = \int_{\partial D_R} (A w)^\alpha
(\overline{Aw})^\beta \partial|Aw|^2 \w ( \ddbar |A
w|^2 )^{n-1} = \lambda^{\alpha} \bar \lambda^{\beta} I_{\alpha,\beta,R}.
\]
It is now clear that $I_{\alpha,\beta,R} = 0$ if $\alpha \neq
\beta$. A similar calculation gives
\[
I_{\alpha,\alpha,R} = \int_{\partial D} (R w)^\alpha (R \bar w)^\alpha
\partial|R w|^2 \w ( \ddbar |Rw|^2 )^{n-1} = R^{2(|\alpha|+n)}
I_{\alpha,\alpha,1},  
\]
which shows that $I_{\alpha,\alpha,R}$ is proportional to
$R^{2(|\alpha|+n)}$. 
\end{proof}

\begin{proof} (of Theorem \ref{linea}) 
The setup here is the same as in Section 4 in \cite{MA2}.
We view $E = T^\ast_{1,0}(\Cn)$ as a rank $n$ vector bundle over $\Cn$ and let
$f(\z) = 2 \pi i(\z - z) \cdot e^\ast$, where $z \in \Cn$ is
fixed, be a section of $E^\ast = T_{1,0}(\Cn)$. We can now express the
contraction 
$\delta_{\z-z}$ as $\delta_f$ operating on $E$. Take a section
$\phi(\z,z)$ of $\Lambda^l E$ with holomorphic polynomials of degree
$r$ for coefficients. We want to show that $\phi = \delta_f
\psi$, where $\psi$ is a section of $\Lambda^{l+1} E$ with holomorphic
polynomials of degree $r-1$ for coefficients. 

Let $D \subset \Cn_w$ be the unit ball. We let $h = dw \w e^\ast$, 
then $\delta_{w-\z} h = 2 \pi i (w-\z) \cdot e^\ast = f(w) - f(\z)$, so
that $h$ will work as a Hefer form. Let 
\[
\sigma = \frac{\bar w \cdot e}{2 \pi i(|w|^2 - \bar w \cdot z)},
\]
then $\delta_{f(w)} \sigma = 1$ if $w \neq z$. Set $u = \sum \sigma \w (\dbar
\sigma)^k$, note that $u$ depends holomorphically on $z$. Since $u$ is
integrable, we let $U$ be the trivial extension over $z$. Moreover, we
have $(\delta_{f(w)} - \dbar)U = 1 - R$, where $R_k = 0$ for $k<n$ and
$(\delta_h)_{n} R_{n} = [z]$, where $[z]$ denotes the $(n,n)$-current
point evaluation at $z$. This follows, for example, from
Theorem 2.2 in \cite{MA4}. Let $\phi(\z,z)$ be a section of
$\Lambda^l E$ with $l>0$ which satisfies $\delta_{f(\z)} \phi = 0$. We set 
\[
g(w,\z) = \delta_{f(\z)} \sum \delta_h^{k-1} (U_k \phi(w,z)) + \sum \delta_h^k
(R^k \phi(w,z));
\]
a calculation shows that $(\delta_{w-\z} - \dbar)g = 0$ and
$g_{0,0}(\z,\z) = \phi(\z,z)$. Note that the second sum is actually
zero since $R \w \phi = 0$. 

Now define 
\[
g_1(w,\z) = \chi_D - \dbar \chi_D \w \frac{s}{\nabla_{w-\z}s} = \chi_D - \dbar
\chi_D \w \sum s \w (\dbar s)^k
\]
where 
\[
s = \frac{1}{2 \pi i} \frac{\bar w \cdot dw}{|w|^2 - \bar w \cdot \z}. 
\]
Note that $(\delta_{w-\z} - \dbar)g_1 = 0$ and
$(g_1)_{0,0}(\z,\z) = 1$. Also, $g_1$ depends holomorphically on $\z$.  

We can now take a current $u = u_{1,0} + \ldots + u_{n,n-1}$ that is smooth
outside $\z$ and such that $(\delta_{w-\z} - \dbar) u = 1-[\z]$. We
have $(\delta_{w-\z} - \dbar)(u \w g \w g_1) = g \w 
g_1 - [\z]\w g \w g_1 = g \w g_1 -
\phi(\z,z)[\z]$, so 
that $\dbar (u \w g \w g_1)_{n,n-1} = \phi(\z,z)[\z] - (g \w
g_1)_{n,n}$. Stokes' theorem gives us $\phi(\z,z) = \delta_{f(\z)}
T\phi$ where  
\[
T\phi =  \int_{\partial D_w} \sum_1^{n-l} \frac{\bar w
  \cdot e \w (\delta_h)_{k-1} [ (d \bar w \cdot e)^{k-1} \w \phi] \w
  \partial |w|^2 \w (\frac{i}{2 \pi} \ddbar |w|^2)^{n-k}}
{2 \pi i (1-\bar w \cdot z)^k (1-\bar w \cdot \z)^{n-k+1}}.
\]
This is seen by noting that $\delta_h \sigma \w s = 0$ and that
\[
\sigma \w \dbar \sigma = \sigma \w \frac{d \bar w \cdot e}{|w|^2 - \bar
w \cdot z}.
\]
First, note that $T\phi$ is holomorphic in $\z$ and $z$ and takes
values in $\Lambda^{l+1} E$. To see that the coefficients are
polynomials of degree $r-1$, we write  
\[
\frac{1}{1 - \bar w \cdot \z} = \sum_0^\infty (\bar w \cdot \z)^k
\]
and similarly for $1/(1 - \bar w \cdot z)$. The only place in $T\phi$
where $w$'s occur is in $\phi$, and according to Lemma
\ref{ineffable} we need to match those up with a corresponding number
of $\bar w$'s. These can come either from the factor $\bar w \cdot
e$, or from the geometric expansions. It is clear that if we start with a
term in $\phi(w,z)$ which is of degree $r$ in $w$, after we
integrate we shall get a term which is of combined degree $r-1$ in $\z$
and $z$. It is easy to see that the same thing is also true if we
start with a term which is of combined degree $r$ in $z$ and $w$. 
The case when $\phi(w)$ takes values in the trivial bundle is proved
in a similar way.  
\end{proof}

We will now find Hefer forms for the complex $E'$ by means of
Theorem \ref{linea}. 

\begin{prop} \label{kowalski}
There exist $(k-l,0)$-form-valued mappings $H_k^l: E_k \to E_l$,
where the coefficients are holomorphic polynomials in $\z$ and $z$,
such that $H_k^l = 0$ for $k<l$, $H_l^l = I_{E_l}$, and 
\be \label{fraser} 
\delta_{\z-z} H_k^l = H_{k-1}^l F_k(\z) - F_{l+1}(z) H_k^{l+1}.
\ee
Moreover, the polynomial degree of $(H_k^l)_{\alpha \beta}$ (that is, the
coefficient of $e_{l \alpha} \otimes e_{k \beta}^\ast$) is
$d_{l}^\alpha-d_{k}^\beta-(k-l)$. 
\end{prop}

\begin{proof}
The proof will be by induction, and by application of Theorem
\ref{linea}. For $k-l=1$ we must solve
$\delta_{\z-z} H_{l+1}^l = F_{l+1}(\z) - F_{l+1}(z)$. We can apply
Theorem \ref{linea} to $(H_{l+1}^l)_{\alpha \beta}$, and it follows
that $(H_{l+1}^l)_{\alpha \beta}$ is a  
$(1,0)$-form whose coefficients are 
holomorphic polynomials of degree $d_{l}^\alpha-d_{l+1}^\beta-1$. 

Assume that the proposition holds for $H_k^l$ with $k-l=i$. We first prove 
that the right hand side of (\ref{fraser}) is $\delta_{\z-z}$-closed. We
have 
\[
\delta_{\z-z} (H_{k-1}^l F_k(\z) - F_{l+1}(z) H_k^{l+1}) = 
\]
\[
= (H_{k-2}^l F_{k-1}(\z) - F_{l+1}(z) H_{k-1}^{l+1})F_k(\z) - F_{l+1}(z)
(H_{k-1}^{l+1} F_k(\z) - F_{l+2}(z) H_k^{l+2}) = 0,
 \]
since $E'$ is a complex. Now regard the   
\[
\delta_{\z-z} (H_k^l)_{\alpha \beta} = (H_{k-1}^l F_k(\z) - F_{l+1}(z)
H_k^{l+1})_{\alpha \beta} = 
\]
\[
= \sum_{\gamma = 1}^{r_{k-1}} (H_{k-1}^l)_{\alpha \gamma}
(F_k(\z))_{\gamma \beta} - \sum_{\gamma = 1}^{r_{l+1}} (F_{l+1}(z))_{\alpha \gamma}
(H_k^{l+1})_{\gamma \beta}.
\]
where $k-l=i+1$. The right hand side is a form of bidegree
$(k-l-1,0)$, and it is easy to see by examining the sums that the
coefficients are 
polynomials of degree $d_{l}^\alpha-d_{k}^\beta - (k-l-1)$. It follows
from Theorem \ref{linea} that $(H_k^l)_{\alpha \beta}$ exists and is a form of
bidegree $(k-l,0)$ whose coefficients are holomorphic polynomials of degree
$d_{l}^\alpha-d_{k}^\beta - (k-l)$.
\end{proof}

We will now prove Proposition \ref{yak} by adapting the Hefer forms
$H_k^l$ to the complex $E$ on $\Pn$. From now on, we change notation
so that $\z = (\z_0, \z') \in \C^{n+1}$ and $z = (1,z') \in
\C^{n+1}$. 

\begin{proof} (of Proposition \ref{yak}) 
The form $(H_k^l)_{\alpha \beta}$  has bidegree $(k-l,0)$, so we have 
\[
(H_k^l)_{\alpha \beta} = \sum_{|I| = k-l} a_I d(\z')_I.
\]
We set 
\[
(h_k^l)_{\alpha \beta} (\z,z') = \z_0^{d_{l}^\alpha-d_{k}^\beta+ (k-l)} \sum_{|I| =
  k-l}a_I (\z'/\z_0,z' ) d ( \z'/\z_0)_I. 
\]
It is clear that $(h_k^l)_{\alpha \beta}$ is a projective
form, since we are multiplying with a high enough degree of $\z_0$,
and that it takes values in the right line bundle. Now, note that 
\[
\delta_\eta d (\z_j/\z_0) = z_j/\z_0 - \z_j/\z_0^2. 
\]
We have  
\[
\delta_\eta (h_k^l)_{\alpha \beta} = \z_0^{d_{l}^\alpha-d_{k}^\beta + (k-l-1)} \sum_{|I| =
  k-l} \sum_{j \in I} \pm a_I (\z'/\z_0,z') (z_j - \z_j/\z_0
) d (\z'/\z_0)_{I \setminus j} = 
\]
\[
= \z_0^{d_{l}^\alpha-d_{k}^\beta + (k-l-1)} \delta_{\z'/\z_0-z'}
(H_k^l)_{\alpha \beta} (\z'/\z_0,z') = 
\]
\[
= \z_0^{d_{l}^\alpha-d_{k}^\beta + (k-l-1)} (H_{k-1}^l (\z'/\z_0,z')
f_k(\z'/\z_0)- f_{l+1}(z') H_k^{l+1} (\z'/\z_0,z'))_{\alpha \beta} = 
\]
\[
= (h_{k-1}^l f_k- f_{l+1}^z h_k^{l+1})_{\alpha \beta},
\]
so we are finished. 
\end{proof}

\section{The Buchsbaum-Rim complex} \label{xavier}

We will now introduce the Buchsbaum-Rim complex, which is a special
case of our complex $E$ over $\Pn$. For more details about this
complex see e.~g.~\cite{MA3}. Let $P: \Cn \to \textrm{Hom}
(\mathbb{C}^m, \mathbb{C}^r)$ be a generically surjective polynomial
mapping and $Z$ the set where $P$ is not surjective. The mapping $P$
will have the role of $F_1$ in the complex (\ref{salieri}). Assume that
$\deg(P^j) \leq d_j$, where $P^j$ is number $j$ in $P$,
and $d_1 \geq d_2 \geq \cdots$. 
Let $p(\z)$ be the matrix with columns $p^j(\z) = \z_0^{d_j}
P^j(\z'/\z_0)$. The mapping $p$
will have the role of $f_1$ in (\ref{salieri}). If
$L_{1}, \ldots, L_m$ are trivial line 
bundles over $\Pn$, with frames $e_1, \ldots,
e_m$, we can define the rank $m$ bundle 
\[
E_1 = (\mathcal{O}(-d_1) \otimes L_1) \oplus \cdots \oplus (\mathcal{O}(-d_n) 
\otimes L_m).
\]
Let $E_0$ be a trivial vector bundle of rank $r$, with the frame
$\{\epsilon_j\}$, and let $p_{ij}$ be the $i$:th element in the column
$p^j$. We can then view 
\[
p = \sum_{i=1}^m \sum_{j=1}^r p_{ij} \epsilon_i \otimes e_j^\ast
\]
as a mapping from $E_1$ to $E_0$.
Contraction with $p$ acts as a mapping $\delta_{p} : E_1 \otimes E_0^\ast
\to \C$. 

Now, note that we can also write $p = \sum p_i \otimes
\epsilon_i$ where $p_i = \sum_j p_{ij} e_j^\ast$ is a section of
$E_1^\ast$. Let 
\[
\det p = p_1 \w \ldots \w p_r \otimes \epsilon_1 \w
\ldots \w \epsilon_r. 
\]
We now get a complex
\[
0 \to E_{m-r + 1} \stackrel{\delta_{p}}{\longrightarrow} \cdots
\stackrel{\delta_{p}}{\longrightarrow} E_3
\stackrel{\delta_{p}}{\longrightarrow} E_2  
\stackrel{\det {p}}{\longrightarrow} E_1
\stackrel{{p}}{\longrightarrow} E_0 \to 0 
\]
over $\Pn$, where for $k \geq 2$ we have 
\[
E_k = \Lambda^{k+r-1} E_1 \otimes S^{k-2} E_0^\ast \otimes \det E_0^\ast,
\] 
and $S^l E_0^\ast$ is the subbundle of $\bigotimes E_0^\ast$
consisting of symmetric $l$-tensors of $E_0^\ast$. 
By $f: E \to E$ we mean the mapping which is either $p$, $\det p$ or
$\delta_{p}$, depending on which $E_k$ we restrict it to. We assume
that $r$ is odd, so that $\dbar \circ \det f = - \det f \circ \dbar$ and
consequently $\dbar \circ f = - f \circ \dbar$. If $r$ is even, one has to
insert changes of sign at some places. 

Next, we can construct currents $U_k$ with bidegree $(0,k-1)$ and 
$R_k$ with bidegree $(0,k)$, with values in $\textrm{Hom}
(E_0,E_k)$, such that $(f - \dbar)U = I_{E_0} -R$, or in other words 
\[
p U_1 = I_{E_0}, \quad (\det p) U_2 - \dbar U_1 = -R_1, \quad
\textrm{and} \quad 
\delta_{p} U_{k+1} - \dbar U_k = -R_k
\]
for $k \geq 2$. See \cite{MA3} for details of the construction. 

Now, recall that we want to solve the division problem $\Phi = P Q$,
where $\Phi = (\Phi_1, \ldots, \Phi_r)$ is an $r$-column of
polynomials of degree $\rho$. The homogenized version of this is $\phi = pq$,
where $\phi$ takes values in $\mathcal{O}(\rho) \otimes E_0$, so that we can
apply Theorem \ref{inuvik}. To do that, we have to
determine the $d_k^j$ of the Buchsbaum-Rim complex. 

For $E_0$ we have
$r_0 = r$ and $d_0^j = 0$ for all $j$, and 
for $E_1$ we have $r_1 = m$ and $d_1^j = -d_j$. Furthermore, for $E_k$
with $k \geq 2$ 
we have $r_k = {m \choose k+r-1} r^{k-2} / (k-2)!$, where ${m \choose k+r-1}$ is
the dimension of $\Lambda^{k+r-1} E_1$ and $r^{k-2} / (k-2)!$ is the dimension of
$S^{k-2} E_0^\ast$. However, it is only the factor $\Lambda^{k+r-1}
E_1$ that contributes with any non-trivial line bundles, and and so we
will not be interested in the factors in $S^{k-2} E_0^\ast \otimes \det
E_0^\ast$. We then have $d_k^{I} = -\sum_{i \in I} d_i$ where $|I|=k +
r -1$, by which we mean that the coefficient of the basis vector $e_I$
in $\Lambda^{k+r-1} E_1$ will take values in $\mathcal{O}(d_k^I) = \mathcal{O}(-\sum_{i
  \in I} d_i)$. 

Recall that $f_k = \delta_p: E_k \to E_{k-1}$ for $k \geq 2$. Let
$f_k^{IJ}$ be the 
coefficent of $e_I \otimes e_J^\ast$, where $|I| = k+r-2$, $|J| = k+r-1$ and
$I \subset J$. Then $f_k^{IJ}$ takes values in 
\[
\mathcal{O}(d_{k-l}^I - d_k^J) = \mathcal{O}(\sum_{i \in I}
d_i -\sum_{i \in J} d_i) = \mathcal{O}(d_j)
\]
where $\{j\} = J \setminus I$. In the same way, $H_k^l: E_k \to
E_{l}$, and for $k,l \geq 2$, we have that $(H_k^l)^{IJ}$ is a $(k-l)$-form
whose coefficients take values in $\mathcal{O}(\sum_{i \in J \setminus I} d_i +
k-l)$ and are holomorphic polynomials in $z$ of degree 
$\sum_{i \in J \setminus I} d_i - (k-l)$. Also, $(H_k^l)^{IJ} = 0$ if
$I \not \subset J$.

\begin{thm} \label{eliot}
Let $P(\z'): \Cn \to \textrm{Hom}
(\mathbb{C}^m, \mathbb{C}^r)$ be a generically surjective polynomial
mapping where the columns $P^j$ have degree $d_j$, and $\Phi(\z') =
(\Phi_1, \ldots, \Phi_r)$ a tuple of polynomials with degree
$\rho$. We have the following decomposition:
\[
\Phi(z') = P(z')Q(z') + 
\int_\z H^0 R \phi \w \alpha^{n + \rho},
\]
where $Q(z') = \sum_1^{m} Q_i(z') e_{j}$ with 
\[
Q_i(z') = \int_\z \z_0^{d_i} (h^1 U \phi)_i \w \alpha^{n + \rho}.
\]
If $R \phi = 0$, the second integral will vanish and we
get a solution to the division problem. Moreover, the polynomial
degree of $PQ$ is $\rho + d_1 + \ldots + d_{\mu + r}$, where
$\mu = \min(n,m-r)$. 
\end{thm}

\begin{proof}
From Theorem \ref{inuvik}, we have 
\[
\deg P Q \leq \max_{k,I} (\rho - d_k^I). 
\]
Since $d_k^{I} = -\sum_{i \in I} d_i$ with $|I| = k+r-1$,
the degree is obviously largest when $k$ is maximal and when we choose
$I = (d_1, \ldots, d_{k+r-1})$. We recall that $U_k = 0$ if $k
> \mu+1$, where $\mu = \min(n,m-r)$, and the theorem follows.
\end{proof}

\begin{preex}
\emph{If $r=1$, the Buchsbaum-Rim complex will simply be the Koszul
complex. To obtain Hefer forms, we first find a
$(1,0)$-form $h(\z,z') = \sum h_{ij} \epsilon_i \otimes e_j^\ast$
taking values in $\mathcal{O}(1) \otimes \Hom (E_1,E_0)$, such that $\delta_\eta 
h = p - p^z$. It is then easy to show that $h_k^l = (\delta_h)^{k-l}$
will act as a Hefer form. We get the decomposition }
\[
\Phi(z') = P(z') Q(z')
+ \int_{\Pn} \sum_k \delta_h^k (R^k \phi(\z)) \w \alpha^{n+\rho},
\]
\emph{where $Q(z') = \sum_1^{m} Q_i(z') e_{j}$ with 
\[
Q_i(z') = \int_\z \z_0^{d_i} \sum_k \delta_h^{k-1} U_k \phi(\z) \w
\alpha^{n + \rho}. 
\]
Furthermore, $\deg P Q \leq \rho + d_1 + \cdots + d_{\mu+1}$, where $\mu =
\textrm{min}(m-1,n)$.}
\end{preex}

\begin{preremark} \label{martha}
\emph{Following Hergoualch \cite{HE}, who uses an idea of Fuhrmann
\cite{FU}, we can reduce the case $r > 1$ to the 
case $r=1$, but we will then lose precision in the degree
estimates. If we wish to solve $P Q = \Phi$, as before, we can
instead set $p' = p_1 \w \ldots \w p_r$ and solve $p'
\cdot q_j' = \phi_j$. Note that we will then get a solution $q_j'$ taking
values in $E_1' = \Lambda^r E_1$. Since the vector space $E_1'$ has
dimension $m' = {m \choose r}$, the previous example will give $\deg
P' Q_j' \leq \rho + d_1' + \cdots + d_{\mu'+1}'$, where $\mu' =
\textrm{min}({m \choose r}-1,n)$. We have $d_1' = d_1 + \cdots + d_r$,
so we can make the estimate 
$\deg P' Q_j' \leq \rho + (\mu'+1)(d_1 + \cdots + d_r)$.}

\emph{Now, we can set 
\[
q_j = (-1)^j \delta_{p_r} \cdots \delta_{p_{j+1}}
\delta_{p_{j-1}} \cdots \delta_{p_1} q_j'.
\]
Note that $q_j$ takes values in $E_1$ and that $p_j q_j =
\phi_j$. Let $q = q_1 + \cdots q_r$, then we have $P Q = \Phi$, with
$\deg P Q \leq \rho + (\mu'+1)(d_1 + \cdots + d_r)$.}
\end{preremark}

\begin{preremark}
\emph{It is possible to get a slightly better result than Theorem~\ref{eliot} 
by means of solving successive $\dbar$-equations in the Buchsbaum-Rim
complex, as is done in \cite{MA3}. We will sketch this approach. 
First, we modify the Buchsbaum-Rim complex by taking tensor products
with $\mathcal{O}(\rho)$:
\[
0 \to E_{m-r + 1} \otimes \mathcal{O}(\rho) \stackrel{\delta_{p}}{\longrightarrow} \cdots
\stackrel{\det p}{\longrightarrow} E_1 \otimes
\mathcal{O}(\rho)\stackrel{p}{\longrightarrow} E_0 \otimes \mathcal{O}(\rho) \to 0,
\]
If $\phi$ is a section of $E_0 \otimes \mathcal{O}(\rho)$ and $R^f \phi = 0$, it
is clear that $v = \phi U$ solves $\nabla_{p} v = \phi$. As in the
classical Koszul complex method of solving division problems, we need
to solve a succession of equations $\dbar w_k = v_k + \delta_{p}
w_{k+1}$. The holomorphic solution will then be $v_1 + \delta_{p} w_2$. 
In order to solve the $\dbar$-equations, we need to assume that 
\[
m \leq n + r - 1 \ \ \ \textrm{or} \ \ \ \rho \geq \sum_1^{n+r} d_j - n, 
\]
where $\deg \Phi
\leq \rho$. With this method, one proves that there exist polynomials
$Q_1, \ldots, Q_m$ such that $\deg P Q \leq \rho$ (Theorem 1.8 in
\cite{MA3}). }

\emph{Actually, by Proposition 5.1 in \cite{EG} we can solve
  $\dbar$-equations for $(0,q)$-forms taking values in $\mathcal{O}(l)$
  explicitly. It is thus possible to obtain explicit  
solutions with these better degree estimates as well, although these
solutions will be   
more complicated than the ones from Theorem~\ref{eliot}, since they
will involve $\mu$ integrations. We do not know if it is possible to
find a weight such that the method in the present article yields
explicit solutions $Q_1, \ldots, Q_m$ such that $\deg P Q \leq \rho$.}
\end{preremark}

\begin{preremark} 
\emph{In the case of the Buchsbaum-Rim complex we can at
  least partially find
explicit Hefer forms. As in the example above, we find $h$ such that
$\delta_\eta h = p - p^z$, 
and we let $h_k^l = \delta_h^{k-l}$ for $k > l \geq 2$. }
\emph{
We also have $h_1^0 = h$ and }
\[
h_2^1 = \sum_{k=1}^r (-1)^{k+1} p^1 \w \ldots \w p^{k-1} \w h^k \w
(p^{k+1})^z \w \ldots \w (p^r)^z \otimes \epsilon_1 \w \ldots \epsilon_r,
\]
\emph{where $h^k$ is the $k$:th column of $h$. We refrain from trying
  to find explicit expressions for $h_k^0$ 
for $k \geq 1$ and $h_k^1$ for $k \geq 2$.} 
\end{preremark}

We can get several corollaries to Theorem \ref{eliot}, that give
conditions which ensure that $R^f \phi = 0$, and state what the estimate of the
degree is in that case. First, there is a classical theorem by
Macaulay \cite{MAC} which says that if $r = 1$, $\Phi = 1$ and $p^1,
\ldots, p^m$ lack common zeroes in $\Pn$, then we have the estimate
$\deg P^i Q_i \leq \sum_1^{n+1} d_j -n$. By means of residue currents,
one can get a generalization to the case $r>1$: 

\begin{prop} 
Assume that $Z$ is empty. Then there exists a matrix $Q$
of polynomials 
such that $PQ  = I_r$, and $\deg P Q \leq \sum_1^{n+r} d_j - n$. 
\end{prop}
\noindent
This is Corollary 1.9 in \cite{MA3}. By using integral
representations, we can get instead

\begin{cor} \label{mountie}
Assume that $Z$ is empty. Then Theorem \ref{eliot} yields an explicit matrix $Q$
of polynomials 
such that $PQ  = I_r$, and $\deg P Q \leq \sum_1^{n+r} d_j$. 
\end{cor}

\begin{proof}
We first note that necessarily $\mu = n$ since $\textrm{codim} Z =
m-r+1$. The result now 
easily follows by several application of Theorem \ref{eliot} where one
chooses the columns of $I_r$ as $\Phi$. 
\end{proof}

If $Z$ is not empty, we would like some condition which guarantees
that $R^f \phi = 0$. We define the pointwise norm $\| \phi \|^2 = \det
(pp^\ast) |p^\ast (pp^\ast)^{-1} \phi |^2$. 

\begin{cor}
If 
\be \label{bryhnia}
\| \phi \| \leq C \det (p p^\ast)^{\min(n,m-r+1)},
\ee
then Theorem \ref{eliot} yields explicit $Q$ such that $PQ = \Phi$ with $\deg P
Q \leq \rho + d_1 + \ldots + d_{\mu + r}$, where $\mu = \min(n,m-r)$.
\end{cor}

This follows from Proposition 1.3 in \cite{MA3}, which says that if
(\ref{bryhnia}) holds, then $R^f \phi = 0$. Just as for the previous
corollary, this is not the optimal result, which states that there
exists a solution $Q$ such that $\deg P Q \leq \rho$ (see Corollary
1.10 in \cite{MA3}). 

\begin{preremark}
\emph{Another possible application for Theorem \ref{inuvik} is the
Eagon-Northcott complex. This is used if we want to solve the equation
$(\det p) \cdot q = \phi$, where $\phi$ is a scalar function. The complex is
given by 
\[
\cdots\stackrel{\delta_{p}}{\longrightarrow} E_3
\stackrel{\delta_{p}}{\longrightarrow} E_2
\stackrel{\delta_{p}}{\longrightarrow} E_1
\stackrel{\det {p}}{\longrightarrow} \C \to 0 
\]
where
\[
E_k = \Lambda^{k+r-1} E_1 \otimes S^{k-1} E_0^\ast \otimes \det E_0^\ast.
\] 
We have already solved this equation in Remark \ref{martha}, but by applying
Theorem \ref{inuvik} to the Eagon-Northcott 
complex we can obtain a sharper degree estimate. Note that in this
division problem we actually have $r=1$ and so we could solve it by
means of the Kozsul complex, but we can improve the degree
estimate by choosing a complex that takes advantage of the fact that
our mapping $\det p$ is constructed from another mapping $p$.}
\end{preremark}


\begin{thebibliography}{10}

\bibitem{MA1} M. Andersson, \emph{The membership problem for
    polynomial ideals in terms of residue currents.}
  Ann. Inst. Fourier (Grenoble)  56  (2006),  no. 1, 101--119. 

\bibitem{MA2} M. Andersson, \emph{Integral representation with
    weights II.} Math. Z. (2006) 254: 315--332.

\bibitem{MA3} M. Andersson, \emph{Residue currents of holomorphic
    morphisms. }  J. Reine Angew. Math.  596  (2006), 215--234.

\bibitem{MA4} M. Andersson, \emph{Integral representation with weights
I.} Math. Ann. 326 (2003), no. 1, pp. 1--18. 

\bibitem{AW} M. Andersson and E. Wulcan, \emph{Residue currents with
    prescribed annihilator ideals.}, Ann. Sci. \'Ecole Norm. Sup. (to
  appear)

\bibitem{BR} Brownawell, W. D., \emph{Bounds for the degrees in the
    Nullstellensatz}, Ann. of Math.(2) 126 (1987),
no. 3, pp. 577--591.

\bibitem{FU} P. A. Fuhrmann, \emph{On the corona theorem and its
    application to spectral problems in Hilbert space.}
  Trans. Amer. Math. Soc.  132  1968 55--66. 

\bibitem{EG} E. Götmark, \emph{Weighted integral formulas on
    manifolds},  to appear in Ark. Mat.

\bibitem{LIC} E. Götmark, \emph{Some Applications of Weighted Integral 
Formulas.}, licentiate thesis, preprint no. 2005:15 from the Department
of Mathematical Sciences, Göteborg University. 

\bibitem{HE} Hergoualch, J., \emph{Le problème de la couronne},
  Mémoire, Bordeaux 2001, Autour du problème de la couronne, Thesis,
  Bordeaux 2004.

\bibitem{KO} J. Kollár, \emph{Sharp Effective Nullstellensatz},
  J. Amer. Math. Soc. 1 (1998), no. 4, pp. 963--975. 

\bibitem{MAC} F. S. Macaulay, \emph{The algebraic theory of modular
    systems}, Cambridge Univ. Press, Cambridge 1916.

\end{thebibliography}
\end{document}